\documentclass[12pt]{amsart}
\PassOptionsToPackage{usenames,dvipsnames}{xcolor}
\usepackage{txfonts}
\PassOptionsToPackage{hyphens}{url}\usepackage{hyperref}
\usepackage{subfig}
\usepackage{amscd,amssymb,mathabx}
\usepackage[arrow,matrix,graph,frame,poly,arc,tips]{xy}
\usepackage{graphicx,calrsfs}
\usepackage{mathtools}
\usepackage{tikz}
\usetikzlibrary{decorations.markings}
\usetikzlibrary{shapes}
\usetikzlibrary{arrows}
\usetikzlibrary{backgrounds}
\usetikzlibrary{calc}
\usepackage{mdwlist}
\usepackage{multirow}
\usepackage{enumerate}
\usepackage[shortlabels]{enumitem}

\usepackage{verbatim}
\usepackage{etoolbox}
\AtEndEnvironment{proof}{\setcounter{claim}{0}}

\DeclarePairedDelimiter{\form}{\langle}{\rangle}

\newcommand\ba{\begin{align*}}
\newcommand\ea{\end{align*}}
\newcommand\be{\begin{enumerate}}
\newcommand\ee{\end{enumerate}}
\newcommand\bp{\begin{proof}}
\newcommand\ep{\end{proof}}
\newcommand\bpp{\begin{prop}}
\newcommand\epp{\end{prop}}
\newcommand\bpb{\begin{prob}}
\newcommand\epb{\end{prob}}
\newcommand\bd{\begin{defn}}
\newcommand\ed{\end{defn}}
\newcommand\bh{\begin{hint}}
\newcommand\eh{\end{hint}}

\newcommand\bN{\mathbb{N}}

\newcommand\bR{\mathbb{R}}

\newcommand\bQ{\mathbb{Q}}

\newcommand\bZ{\mathbb{Z}}

\newcommand\supp{\operatorname{supp}}

\newcommand\Mod{\operatorname{Mod}}

\DeclareMathOperator\Homeo{Homeo}

\DeclareMathOperator\CritReg{CritReg}

\newcommand\sse{\subseteq}
\newcommand\co{\colon}

\DeclareMathOperator\Fix{Fix}
\DeclareMathOperator\Out{Out}
\DeclareMathOperator\Aut{Aut}

\DeclareMathOperator\Diff{Diff}

\newcommand\rot{\operatorname{rot}}

\renewcommand{\MR}[1]
{\href{http://www.ams.org/mathscinet-getitem?mr=#1}{MR#1}}

\def\thetitle{{Virtual critical regularity of mapping class group actions on the circle}}
\def\theauthors{{Sang-hyun Kim, Thomas Koberda, Crist\'obal Rivas}}
\usepackage{hyperref}
\hypersetup{
   colorlinks=false,
   plainpages,
   urlcolor=black,
   linkcolor=black
   pdftitle=   \thetitle,
   pdfauthor=  {\theauthors}
}

\theoremstyle{plain}
\newtheorem{thm}{Theorem}[section]
\newtheorem{lem}[thm]{Lemma}
\newtheorem{cor}[thm]{Corollary}
\newtheorem{prop}[thm]{Proposition}

\newtheorem{que}[thm]{Question}
\newtheorem*{principle*}{Principle}
\newtheorem*{claim*}{Claim}

\theoremstyle{remark}

\newtheorem{rem}[thm]{Remark}

\theoremstyle{definition}
\newtheorem{defn}[thm]{Definition}
\newtheorem{prob}{Problem}[section]

\begin{document}
\title\thetitle
\date{\today}

\author[S. Kim]{Sang-hyun Kim}
\address{School of Mathematics, Korea Institute for Advanced Study (KIAS), Seoul, 02455, Korea}
\email{skim.math@gmail.com}
\urladdr{https://www.cayley.kr}

\author[T. Koberda]{Thomas Koberda}
\address{Department of Mathematics, University of Virginia, Charlottesville, VA 22904-4137, USA}
\email{thomas.koberda@gmail.com}
\urladdr{http://faculty.virginia.edu/Koberda}

\author[C. Rivas]{Crist\'obal Rivas}
\address{Dpto. de Matem\'aticas y C.C., Universidad de Santiago de Chile, Alameda 3363, Santiago, Chile}
\email{cristobal.rivas@usach.cl}
\urladdr{{http://mat.usach.cl/index.php/2012-12-19-12-50-19/academicos}\linebreak/183-cristobal-rivas}

\begin{abstract}
We show that if $G_1$ and $G_2$ are non-solvable groups, then  no $C^{1,\tau}$ action of $(G_1\times G_2)*\bZ$ on $S^1$ is faithful for $\tau>0$.
As a corollary, if $S$ is an orientable surface of complexity at least three then the critical regularity of an arbitrary finite index subgroup of the
mapping class group $\Mod(S)$ with respect to the circle is at most one, thus strengthening a result of the first two authors with Baik.
\end{abstract}

\maketitle
\section{Introduction}

Let $G$ be a group, and let $M$ be a smooth manifold. For $k\in\bN$ and $\tau\in[0,1]$, we denote by 
$\Diff^{k,\tau}(M)_0$ the group of $C^k$ diffeomorphisms of $M$ whose $k^{th}$ derivatives are
$\tau$--H\"older continuous and are isotopic to the identity.

The \emph{critical regularity of $G$ with respect to $M$} is defined to be
\[\CritReg_M(G)=\sup\{k+\tau\mid k\in\bN,\tau\in[0,1]\text{ and }G \textrm{ injects into } \Diff^{k,\tau}(M)_0\}.\]
 By convention,
$\Homeo(M)_0=\Diff^0(M)_0$, and if $G$ admits no injective homomorphism into $\Homeo(M)_0$ then $\CritReg_M(G)=-\infty$.

\subsection{Main results}
In this article, we concentrate on computing the critical regularity of certain groups in the case $M=S^1$, and we will
suppress $M$ from the notation; therefore, we write
\[\CritReg(G):=\CritReg_{S^1}(G).\]
Note that $\Homeo(S^1)_0=\Homeo_+(S^1)$, where the right hand side denotes the group of 
orientation preserving homeomorphisms of $S^1$. Our main result
is as follows.

\begin{thm}\label{thm:main}
If $G_1$ and $G_2$ are non-solvable groups,
then 
\[\CritReg((G_1\times G_2)*\bZ)\leq 1.\]
\end{thm}

Every countable subgroup $G$ of $\Homeo^+(S^1)$ is topologically conjugate to a group of 
bi--Lipschitz homeomorphisms of $S^1$ by~\cite{DKN2007}.
Moreover, the group $G\ast \bZ$ admits an embedding into $\Homeo^+(S^1)$ by~\cite{BS2018}.
It follows that
\[
\CritReg(G)=\CritReg(G*\bZ)\ge 1.\] 
The following is now an immediate corollary of the main theorem.

\begin{cor}\label{cor:f2xf2}
We have \[\CritReg((F_2\times F_2)*\bZ)= 1.\]
\end{cor}

We note that the group $(F_2\times F_2)\ast \bZ$ admits a faithful $C^1$--action on $S^1$ and on $I:=[0,1]$,
as does every finitely generated residually torsion--free nilpotent group~\cite{FF2003,Jorquera,Parkhe2016}, and so Corollary~\ref{cor:f2xf2} is
optimal.

Corollary~\ref{cor:f2xf2} allows us to compute the critical regularity of many mapping class groups of surfaces. Recall that if $S$ is an
orientable surface of genus $g$ and with $n$ punctures, boundary components, and marked points, we write $\Mod(S)$ for the group
of isotopy classes of homeomorphisms of $S$ that preserve the punctures, boundary components, and marked points (pointwise). We use
$\xi(S)$ for the \emph{complexity} of $S$, which is defined by \[\xi(S)=3g-3+ n.\]

If $g\geq 2$
and $n=1$ then $\Mod(S)$ acts faithfully on $S^1$, and if $S$ has a boundary component then $\Mod(S)$ acts faithfully on 
$I$~\cite{CassonBleiler,Nielsen1927,HT1985}. It was shown in~\cite{FF2020} that the critical regularity of $\Mod(S)$ is at most two, provided
that $g\geq 3$. This was strengthened in~\cite{Parwani2008,MannWolff20}, where it was shown that the critical regularity of $\Mod(S)$ is
at most one. These latter results in fact showed that any $C^1$ action of the full mapping class group on
$S^1$ factors through a finite group.

For finite index subgroups $H<\Mod(S)$, the critical regularity question is more complicated
because finite index subgroups of mapping class groups are poorly understood.
The first two authors and Baik~\cite{BKK2019JEMS} proved that if $\xi(S)\geq 2$,
then every finite index subgroup $H$ of $\Mod(S)$
satisfies $\CritReg(H)\leq 2$, answering a question of Farb in~\cite{Farb2006}. 
In~\cite{Parwani2008}, it is shown that \emph{if} every finite index subgroup
of the mapping class group has finite abelianization when $g\geq 3$ (i.e.~\emph{if} the Ivanov Conjecture holds), then $\CritReg(H)\leq 1$ for
$H<\Mod(S)$ of finite index and $S$ of genus at least $6$ (and in fact no faithful $C^1$ action exists).

Whereas Corollary~\ref{cor:f2xf2} does not rule out the existence of a faithful $C^1$ action of a finite index subgroup of the mapping class
group, it does show that the critical regularity of a finite index subgroup of $\Mod(S)$ is bounded above by one.

\begin{cor}\label{cor:mcg}
Let $S$ be a surface with $\xi(S)\ge3$, and let $\tau>0$.
If $H$ is a finite index subgroup of $\Mod(S)$ then
it admits no faithful $C^{1,\tau}$--action on the circle; in particular $\CritReg(H)\leq 1$.
\end{cor}

To see how Corollary~\ref{cor:f2xf2} implies Corollary~\ref{cor:mcg}, note that under the hypotheses on $S$, there are two subsurfaces
$S_1$ and $S_2$ of $S$ which are homeomorphic to tori with a single boundary component. The mapping class groups of both $S_1$
and $S_2$ contain copies of the group $F_2$,
and the corresponding copies of $F_2$ commute with each other. Adjoining a pseudo-Anosov mapping
class $\psi$ of $S$ gives a $5$--generated group, which after passing to powers of the generators if necessary, furnishes a copy of 
$(F_2\times F_2)*\bZ$ in $\Mod(S)$ (by the main result of~\cite{Koberda2012}, for instance). Corollary~\ref{cor:mcg} then follows
immediately.

Note that Corollary~\ref{cor:f2xf2} also implies an analogous version of Corollary~\ref{cor:mcg} for the groups $\Aut(F_n)$ and
$\Out(F_n)$, since whenever $n\geq 4$,
these groups contain copies of mapping class groups that fall under the purview of Corollary~\ref{cor:mcg}.
We will not comment on these points any further, since unlike mapping class groups of surfaces, automorphism
groups of free groups are not known to have any natural actions on the circle.

We emphasize that Corollary~\ref{cor:mcg} does not rule out the possibility that such a subgroup $H$ may admit a faithful action on
$S^1$ by $C^1$ diffeomorphisms. Whether or not such an action exists appears to be beyond the reach of current technology. The reader
may consult section 8.3.5 of~\cite{KK21-book}, for instance. For certain restricted classes of finite index subgroups, it is known that $H$
cannot act faithfully by diffeomorphisms; see~\cite{parwani2021}.

It is currently an open question for which surface $S$ a finite index subgroup of $\Mod(S)$ admits a faithful $C^0$--action on $S^1$.
In the case when $H$ is such a finite index subgroup we have from the above corollary that $\CritReg(H)=1$. We note that
it is usually quite difficult to compute the critical regularity of a particular group whose critical regularity is known to be finite.
For a survey of results, the reader is directed to~\cite{Navas2008GAFA,CJN2014,JNR2018,KK2020crit,Mann:aa}.

Corollary~\ref{cor:mcg} follows immediately from Corollary~\ref{cor:f2xf2} after observing that under the assumption that
$\xi(S)\geq 3$, the group $\Mod(S)$ and all of its finite index subgroups contain copies of $(F_2\times F_2)*\bZ$ 
(cf.~\cite{Koberda2012,BKK2019JEMS,KKR2020}).

In the case where $\xi(S)\le 1$, the mapping class group of $S$ is virtually free,
so that $\CritReg(H)=\infty$ for a suitable finite index subgroup $H$ of $\Mod(S)$.
The only  case that is left out from Corollary~\ref{cor:mcg} is exactly when $\xi(S)=2$:
\begin{que}
Let $S$ be a twice--punctured torus
or a five--times punctured sphere.
Does some finite index subgroup of $\Mod(S)$ admit a faithful $C^{1,\tau}$ action on $S^1$ with $\tau>0$?
\end{que}

\subsection{A dynamical perspective on the main result}
For the remainder of this section, we frame the discussion of this article in a more precise manner, and while doing so introduce some
relevant concepts. Let  $G$ be a group acting on a space $X$, and
define the \emph{(open) support} of $g\in G$ by
\[
\supp g:=X\setminus\Fix g.\]
The \emph{support of $G$}  is the set 
\[\supp G:=\bigcup_{g\in G}\supp g.\]
We call each point in 
\[\Fix G:=\bigcap_{g\in G}\Fix g\]
a \emph{global fixed point} of $G$.

Let us say $G$ admits a \emph{disjointly supported pair} (or, $G$ is \emph{non--overlapping}) 
if  there exist nontrivial elements $g,h\in G$  satisfying  \[\supp g\cap\supp h=\varnothing.\] 

In~\cite{KKR2020}, the authors proved the following result, which was partially
based on the methods in~\cite{Navas2008GAFA}:
\begin{thm}[\cite{KKR2020}, Theorem 1.1]\label{thm:kkr}
If $G_1$ and $G_2$ are non-solvable groups, 
and if $\tau>0$,
then there is no faithful
$C^{1,\tau}$ action of $(G_1\times G_2)*\bZ$ on a compact interval.
\end{thm}

In that paper, the main technical result was the following.
\begin{thm}[{\cite[Section 4.1]{KKR2020}}]\label{thm:kkr-0}
Let $\tau>0$ be a real number, and let $k\ge3$ be an integer
such that $\tau(1+\tau)^{k-2}\ge1$.
If $G_1$ and $G_2$ are groups that are not solvable of degree at most $k$,
and if
\[G:=G_1\times G_2\longrightarrow\Diff_+^{1,\tau}([0,1])
\]
is an embedding, 
then $G$ contains a disjointly supported pair.
\end{thm}

Theorem~\ref{thm:kkr} follows from Theorem~\ref{thm:kkr-0} by an application of the
$abt$--Lemma from~\cite{KK2018JT}:
\begin{prop}[The $abt$--Lemma]\label{prop:abt}
Let $M\in\{I,S^1\}$ and let $a,b\in\Diff^1_+(M)$ be such that \[\supp a\cap\supp b
=\varnothing.\] Then if $t\in \Diff^1_+(M)$ is arbitrary, the group $\form{ a,b,t}$ is not isomorphic to $\bZ^2\ast \bZ$.
\end{prop}

Proposition~\ref{prop:abt} implies that if a group $G$ {\em always} has  elements with
disjoint supports whenever acting on $I$ or $S^1$ by diffeomorphisms of some regularity, then $G*\bZ$ \emph{never}
acts faithfully by diffeomorphisms on $I$ or $S^1$ of that regularity.
Thus, to prove Theorem~\ref{thm:main}, it will suffice for us to establish the following:

\begin{prop}\label{prop:main}
Let $\tau>0$, and let $G_1$ and $G_2$ be non-solvable groups.
If $\phi$ is a $C^{1,\tau}$--action of $G_1\times G_2$ on $S^1$,
then $G$ admits a disjointly supported pair.
\end{prop}

We will argue Proposition~\ref{prop:main} by showing directly that the commutator subgroup of $G_1\times G_2$
admits a global fixed point, thus reducing to Theorem~\ref{thm:kkr}.

\section{Preliminaries}\label{s:prelim}
For a direct product of groups
\[
G=G_1\times G_2,\]
we identify
$G_1$ with $G_1\times\{1\}$ and $G_2$ with $\{1\}\times G_2$, so that $G_i$ is a normal subgroup of $G$ for $i\in\{1,2\}$;
moreover, we have
$G=\form{ G_1,G_2}$.

Now, suppose a subgroup $G\le\Homeo^+(S^1)$ is given. 
A Borel probability measure $\mu$ on $S^1$ is said to be \emph{$G$--invariant}  if for all $g\in G$ and for all measurable $A\subset S^1$, we have $\mu(A)=\mu(g^{-1}A)$. The \emph{support}
of  $\mu$, denoted as $\supp\mu$, means the largest closed subset $X\subset S^1$ such that every open subset of $X$ has positive measure.

Recall that the \emph{rotation number} \[\rot\colon\Homeo_+(S^1)\longrightarrow\bR/\bZ\] is defined as follows. Let $f\in \Homeo_+(S^1)$,
and lift $f$ to $F\in\Homeo_+(\bR)$. Note that such a lift is always periodic, and that any two such lifts differ by an integer translation. One
chooses an arbitrary $x\in\bR$ and writes \[\rot(f)=\lim_{n\to\infty}\frac{F^n(x)}{n}\pmod{\bZ}.\] It is not difficult to check that the definition
is independent of all the choices made. 

A standard fact is that an orientation preserving
homeomorphism of $S^1$ has nonzero rotation number if and only if it has no fixed points~\cite{athanassopoulos,Navas2011}.
We will appeal to the following basic fact relating rotation numbers and invariant measures; note that the second part of the proposition is an immediate consequence of the first.

\begin{prop}[See~\cite{Navas2011}, Theorem 2.2.10]\label{prop:rot-homo}
If $G\le \Homeo_+(S^1)$ admits an invariant measure $\mu$,
then the restriction
\[
\rot\restriction_G\co G\to \bR/\bZ\]
is a group homomorphism satisfying
\[
\rot(g)=\mu[x,g(x))\]
for all $g\in G$ and $x\in S^1$.
Moreover, the kernel of this homomorphism fixes every point in $\supp\mu$.
\end{prop}

We now recall some ideas from~\cite{KKR2020} that will be crucial in the proof of our main result.
Following \cite{NR2009}, we say that two elements $f,g\in\Homeo_+(\bR)$ are \emph{crossed} if there exist point $u<w<v$ in $\bR$ such that:
\begin{enumerate}
\item
$g^n(u)<w<f^n(v)$ for all $n\in\bZ$;
\item
There is an $N\in\bZ$ such that $g^N(v)<w<f^N(u)$.
\end{enumerate}
A group action of $G$ on $\bR$ by orientation preserving homeomorphisms is called \emph{Conradian} if it admits no crossed elements.

\begin{lem}[\cite{KKR2020}, Lemma 3.10 (Centralizer--Conradian Lemma)]\label{lem:tau-conrad}
Let $\tau>0$, and let $G\le \Diff_+^{1,\tau}([0,1])$. 
If $c$ is a central element of $G$,
then the restriction of $G$ to $\supp c$ is Conradian.
\end{lem}

The relationship between $C^{1,\tau}$ actions and Conradian actions is elucidated by the following technical fact.

\begin{lem}[\cite{KKR2020}, Lemmas 3.4]\label{lem:tau-solv1}
If $\tau,u>0$ are real numbers, and if $k\ge2$ is an integer
satisfying $\tau(1+\tau)^{k-2}\ge u$,
then $\Diff_+^{1,\tau}([0,1])$ does not contain a $(k,u)$--nesting.
\end{lem}

Briefly speaking, a \emph{$(k,u)$--nesting} is a finite set $S\sse\Homeo^+([0,1])$
such that 
for some infinite sequence $(s_1,s_2,\ldots)$ of elements from $S$,
for some nested open intervals
\[ J_1\supsetneq J_2\supsetneq\cdots \supsetneq J_k,\]
and for some choices \[\{t_{n,i}\mid 2\le i\le k, n\ge0\}\sse S,\]
one has 
\[\sum_{n\ge0} |s_n\cdots s_2s_1J_{1}|^u<\infty,\] 
together with
\[
t_{n,i}w_n J_i\cap w_n J_i=\varnothing,\quad
t_{n,i}w_n J_{i-1}=w_n J_{i-1},\] where here $w_n=s_n\cdots s_1$.
A $(k,u)$--nesting is a feature of an action that is weaker than the classical notion of a ``$k$--level structure''~\cite{Navas2008GAFA}.

\begin{lem}[\cite{KKR2020}, Lemma 3.13]\label{lem:tau-solv2}
Let $G\le\Homeo_+([0,1])$ be a Conradian group such that $G^{(k)}\ne1$ for some $k\ge2$.
If $c$ is a central element of $G$ fixing no points in $(0,1)$,
then $G$ contains a $(k,1)$--nesting.
\end{lem}

One may take the $(k,u)$--nesting as a black box for the purpose of this paper, 
and only note the following immediate consequence of the three preceding lemmas.
\begin{lem}\label{lem:tau-solv}
Let $\tau>0$ be a real number and $k\ge3$ be an integer
such that $\tau(1+\tau)^{k-2}\ge1$.
If $c$ is a central element of $G\le\Diff_+^{1,\tau}([0,1])$
fixing no points in $(0,1)$,
then $G^{(k)}=1$.\end{lem}

A fixed point $a$ of $g\in\Diff_+^1(S^1)$ is called a \emph{hyperbolic fixed point}
if $g'(a)\ne1$. The following deep theorem of Deroin--Kleptsyn--Navas (which is a generalization of a result due to Sacksteder)
will be an important ingredient for us.
\begin{thm}[\cite{DKN2007}]\label{thm:dkn}
If a subgroup $G$ of $\Diff_+^1(S^1)$ preserves no probability measure on $S^1$,
then $G$ contains an element $g$ such that $\Fix g$ is nonempty, finite,
and consists entirely of hyperbolic fixed points.
\end{thm}

\begin{rem}\label{rem:hyp-fixed}
For a group $G\le\Homeo^+(S^1)$ that does not admit a finite orbit,
there uniquely exists a smallest, nonempty, closed $G$--invariant set $\Lambda_G$, called the \emph{limit set} of $G$; see Theorem 2.1.1 in~\cite{Navas2011}, for instance.
The limit set is either $S^1$ or a Cantor set, the latter of which is called
 the \emph{exceptional minimal set} of $G$.
In Theorem~\ref{thm:dkn}, we can find a point $x\in \Lambda_G\setminus \Fix g$,
since the limit set $\Lambda_G$ is necessarily infinite.
Consider now the component $J$ of $\supp g$ containing $x$.
The $G$--invariance of $\Lambda_G$ implies that 
\[
\partial J=g^{\pm\infty}(x)\sse\Fix g\cap \Lambda_G.\]
In other words, we can always find a hyperbolic fixed point of $g$ in $\Lambda_G$.
\end{rem}

\section{Establishing the main result}\label{s:main}
A group $G$ is said to be \emph{solvable of degree at most $k$}
if the subgroup $G^{(k)}$, the $k$--th term in the derived series, is trivial.
As noted in the introduction
Theorem~\ref{thm:main} will follow from
Proposition~\ref{prop:main}, which in turn is an immediate consequence of the   stronger result given below.

\begin{thm}\label{thm:main2}
Let $k\ge3$ be an integer,
and let $\tau>0$ be a real number
satisfying $\tau(1+\tau)^{k-2}\ge1$.
If $G_1$ and $G_2$ are groups that are not solvable of degree at most $(k+1)$,
then 
every faithful $C^{1,\tau}$--action of $G_1\times G_2$ on $S^1$ admits a disjointly supported pair.
In particular, we have that
\[
\CritReg\left((G_1\times G_2)\ast\bZ\right)\le 1+\tau.\]
\end{thm}

Note that the second part of the theorem follows from the first along with the $abt$--Lemma (Proposition~\ref{prop:abt}). The lemma below is a key step in the proof of the first part.
\begin{lem}\label{lem:no-meas-solv}
Let $k$ and $\tau$ be as in Theorem~\ref{thm:main2}.
If a group $H\le\Diff_+^{1,\tau}(S^1)$ can be written as a direct product
$H=H_1\times H_2$,
and if $H_1$ does not preserve a probability measure on $S^1$,
then $H_2$ is solvable of degree at most $(k+1)$.
\end{lem}
\bp
From Theorem~\ref{thm:dkn} and Remark~\ref{rem:hyp-fixed}, we can find some $c\in H_1$ 
and $a\in \Lambda_{H_1}\cap\Fix c$ such that 
$c$ fixes finitely many points
and such that $c'(a)\ne1$.
For all $h\in H_2$, the point  $h(a)$ is also a hyperbolic fixed point of $c$ with the derivative $c'(a)$ since
\[
c'\circ h(a)=(c\circ h)'(a)/h'(a)=h'(c(a))\cdot c'(a)/h'(a)=c'(a).\]
It follows that $H_2(a)$ does not have an accumulation point, and in particular is finite.
As $H_2$ admits an invariant probability measure (with atoms at points of $H_2(a)$),
we see from Proposition~\ref{prop:rot-homo} that $K:=[H_2,H_2]$ fixes the point $a$.

Let $U_1$ and $U_2$ be the two components of $\supp c$ containing $a$ on their boundaries.
The group $K$ preserves each $U_i$, since $K$ permutes the components of $\supp c$ and fixes the point $a$.
Applying Lemma~\ref{lem:tau-solv} to the restriction
\[
\left(\form{c}\times K\right)\restriction_{\overline{U_i}},\]
we see that $K^{(k)}$ acts trivially on $U_i$ for $i=1,2$.

Suppose $V$ is a component of the support of $K$, not intersecting $U_1\cup U_2$.
Since $a$ lies in the limit set of $H_1$, we can find some $h_1\in H_1$ such that 
\[h_1(V)\sse U_1\cup U_2.\]
Let $g\in K^{(k)}$ and $v\in V$ be arbitrary.
Since $g$ acts trivially on $h_1(v)$, we have that
\[
g(v) = h_1^{-1}\circ g\circ h_1(v)=h_1^{-1}\circ h_1(v)=v.\]
Combined with the preceding paragraph, this proves that 
\[H_2^{(k+1)}=K^{(k)}=1.\qedhere\]
\ep

We note the following general observation regarding topological actions.
The first part of the lemma is well-known. Presumably, the second part is also known to experts, but the authors could not find a written reference for it.
\begin{lem}\label{lem:c0-product}
Let $H=H_1\times H_2$ be a subgroup of $\Homeo^+(S^1)$.
\be
\item\label{p:fix} If each $H_i$ admits a global fixed point, then so does $H$.
\item\label{p:meas} If each $H_i$ preserves a Borel probability measure on $S^1$, then so does $H$.
\ee
\end{lem}
\bp
(\ref{p:fix})
Suppose not. Since $\Fix H_1\cap \Fix H_2=\varnothing$, 
we can find some $b\in \Fix H_1\cap \supp H_2$.
Let $J$ be the component of $\supp H_2$ containing $b$.
There exists a sequence $\{h_n\}$ in $H_2$ such that 
\[
b':=\lim_{n\to\infty} h_n(b)\in \partial J.\]
Then $b'\in \Fix H_1\cap \Fix H_2$, which is a contradiction.

(\ref{p:meas})
Let $\mu_i$ be a probability measure preserved by $H_i$.
By Proposition~\ref{prop:rot-homo}, 
the restriction of $\rot$ to each $H_i$ is a homomorphism.

Suppose first that $\rot(H_1)\cup \rot(H_2) $
is a  discrete subset of $\bR/\bZ$. This means that $K_i$, the kernel of the map
$\rot: H_i\longrightarrow  \bQ/\bZ$, has finite index in $H_i$ ($i=1,2$). Since each $K_i$ admits a global fixed point,
so does $K_1\times K_2$. This latter group has finite index in $H$,
 and so $H$ has a finite orbit and preserves a probability measure.

We now assume that $\rot(H_1)\cup \rot(H_2) $ is indiscrete in $\bR/\bZ$. Without loss of generality, $\rot(H_1)$
is a dense subgroup of $\bR/\bZ$.  By a result of Plante (See Proposition 2.2 of \cite{plante83}), it follows that $H_1$ preserves a {\em unique} Borel probability measure $\mu_1$. Finally, if $h_2\in H_2$ and $h_1\in H_1$, then
\[h_1^*h_2^*\mu_1=h_2^*h_1^*\mu_1=h_2^*\mu_1.\]
The uniqueness of $\mu_1$ implies that $h_2^*\mu_1=\mu_1$.
In other words we have shown that $\mu_1$ is also $H_2$--invariant, and so also $H$--invariant.  \ep

\bp[Proof of Theorem~\ref{thm:main2}]
We may assume that the given group $G:=G_1\times G_2$ is a subgroup
of $\Diff_+^{1,\tau}(S^1)$.
If some $G_i$ does not admit an invariant probability measure, we apply Lemma~\ref{lem:no-meas-solv} to obtain a contradiction.
So, we will assume that each $G_i$ preserves a probability measure.
Lemma~\ref{lem:c0-product} implies that $G$ also preserves a probability measure $\mu$.

By Proposition~\ref{prop:rot-homo} the rotation number is trivial on the group
\[
H:=[G,G]=[G_1,G_1]\times [G_2,G_2].\]
Moreover, the support of $\mu$ is contained in the global fixed point set of $H$, which is therefore nonempty.
So, the inclusion $H\hookrightarrow \Diff_+^{1,\tau}(S^1)$
factors through an injection $H\hookrightarrow\Diff_+^{1,\tau}([0,1])$.
By Theorem~\ref{thm:kkr-0}, it follows that $H$ admits a disjointly supported pair.
Along with the $abt$--Lemma (Proposition~\ref{prop:abt}), this completes the proof.\ep


\section*{Acknowledgements}
The first author is supported by Samsung Science and Technology Foundation under Project Number SSTF-BA1301-51
and
 by a KIAS Individual Grant (MG073601) at Korea Institute for Advanced Study.
The second author is partially supported  by an Alfred P. Sloan Foundation Research Fellowship, by
NSF Grant DMS-1711488, and by NSF Grant DMS-2002596.
The third author was partially supported by FONDECYT 1210155.

\section*{Statements and Declarations }

The authors declare that they have no conflict of interest.

\bibliographystyle{amsplain}
\providecommand{\bysame}{\leavevmode\hbox to3em{\hrulefill}\thinspace}
\providecommand{\MR}{\relax\ifhmode\unskip\space\fi MR }
\providecommand{\MRhref}[2]{%
  \href{http://www.ams.org/mathscinet-getitem?mr=#1}{#2}
}
\providecommand{\href}[2]{#2}

\end{document}